\documentclass[journal,transmag]{IEEEtran}

\usepackage{graphicx}
\usepackage{xcolor}
\usepackage{hyperref}
\usepackage{amsmath,amssymb}
\usepackage{tikz,circuitikz}
\usetikzlibrary{calc}
\usepackage{pgfplots}
\usepgfplotslibrary{groupplots}
\usepackage{tuda-pgfplots}
\usepackage{pgfplotstable}
\pgfplotsset{compat=1.18}
\usepackage[exponent-product = \cdot]{siunitx}
\usepackage{subcaption}
\usepackage{cite}

\usetikzlibrary{calc,decorations.markings,shapes.misc,positioning,matrix,shapes.geometric}
\usetikzlibrary{patterns,patterns.meta}
\usetikzlibrary {quotes, arrows.meta}
\usetikzlibrary {graphs}
\usetikzlibrary {graphs.standard}
\usetikzlibrary{arrows.meta,matrix}
\usetikzlibrary{positioning}

\newcommand{\CURL}{\ensuremath{\nabla\times}}
\newcommand{\DIV}{\ensuremath{\nabla\cdot}}
\newcommand{\GRAD}{\ensuremath{\nabla}}
\newcommand{\permittivity}{\ensuremath{\epsilon}}
\newcommand{\reluctivity}{\ensuremath{\nu}}
\newcommand{\conductivity}{\ensuremath{\sigma}}

\newcommand{\imNum}{\ensuremath{i}}
\newcommand{\Afield}{\ensuremath{\boldsymbol{A}}} %
\newcommand{\Bfield}{\ensuremath{\boldsymbol{B}}} %
\newcommand{\Dfield}{\ensuremath{\boldsymbol{D}}} %
\newcommand{\Efield}{\ensuremath{\boldsymbol{E}}} %
\newcommand{\Jfield}{\ensuremath{\boldsymbol{J}}} %
\newcommand{\rhofield}{\ensuremath{\varrho}} %
\newcommand{\phifield}{\ensuremath{\varphi}} %

\newcommand{\nfield}{\ensuremath{\boldsymbol{n}}}
\newcommand{\zerofield}{\ensuremath{\boldsymbol{0}}}

\newcommand{\diffV}{\ensuremath{\mathrm{d}\mathrm{V}}}

\newcommand{\scalarWeakSpace}{\ensuremath{H^1(\Omega)}}%
\newcommand{\vectorWeakSpace}{\ensuremath{H(\Omega;\operatorname{curl})}}%

\newcommand{\scalarBaseFun}{\ensuremath{v}}
\newcommand{\vectorBaseFun}{\ensuremath{\boldsymbol{w}}}
\newcommand{\stiffMatCurl}[1]{\ensuremath{\mathbf{K}_{#1}}}
\newcommand{\stiffMatGrad}[1]{\ensuremath{\mathbf{K}}_{#1}}
\newcommand{\massMatCurl}[1]{\ensuremath{\mathbf{M}_{#1}}}
\newcommand{\gradMat}[1]{\ensuremath{\mathbf{G}_{#1}}}
\newcommand{\divMat}[1]{\ensuremath{\mathbf{S}_{#1}}}
\newcommand{\scalarDOF}{\ensuremath{u}}
\newcommand{\scalarDOFs}{\ensuremath{\mathbf{\scalarDOF}}}
\newcommand{\vectorDOF}{\ensuremath{a}}
\newcommand{\vectorDOFs}{\ensuremath{\mathbf{\vectorDOF}}}
\newcommand{\norm}[2]{\left\lVert#1\right\rVert_{#2}}
\newcommand{\scalarDiscreteSpace}{\ensuremath{V}}%
\newcommand{\vectorDiscreteSpace}{\ensuremath{W}}%
\newcommand{\dt}[1]{\ensuremath{\dot{#1}}}%
\newcommand{\dtt}[1]{\ensuremath{\ddot{#1}}}%
\newcommand{\deltaT}{\ensuremath{\Delta t}}%
\newcommand{\LTwoErr}{\ensuremath{\varepsilon_{\mathrm{L2}}}}%

\makeatletter

\def\env@cases#1{%
  \let\@ifnextchar\new@ifnextchar
  \left\lbrace\def\arraystretch{1.2}%
  \array{@{}#1@{\quad}l@{}}}
\makeatother

\tikzset{
	cuboid/.pic={
		\tikzset{%
			every edge quotes/.append style={midway, auto},
			/cuboid/.cd,
			#1
		}
		\draw [every edge/.append style={pic actions, densely dashed, opacity=.5}, pic actions]
		(0,0,0) coordinate (o) -- ++(-\cubescale*\cubex,0,0) coordinate (a) -- ++(0,-\cubescale*\cubey,0) coordinate (b) edge coordinate [pos=1] (g) ++(0,0,-\cubescale*\cubez)  -- ++(\cubescale*\cubex,0,0) coordinate (c) -- cycle
		(o) -- ++(0,0,-\cubescale*\cubez) coordinate (d) -- ++(0,-\cubescale*\cubey,0) coordinate (e) edge (g) -- (c) -- cycle
		(o) -- (a) -- ++(0,0,-\cubescale*\cubez) coordinate (f) edge (g) -- (d) -- cycle;
	},
	/cuboid/.search also={/tikz},
	/cuboid/.cd,
	width/.store in=\cubex,
	height/.store in=\cubey,
	depth/.store in=\cubez,
	units/.store in=\cubeunits,
	scale/.store in=\cubescale,
	width=10,
	height=10,
	depth=10,
	units=cm,
	scale=.1,
}

\tikzset{
	annotated cuboid/.pic={
		\tikzset{%
			every edge quotes/.append style={midway, auto},
			/cuboid/.cd,
			#1
		}
		\draw [every edge/.append style={pic actions, densely dashed, opacity=.5}, pic actions]
		(0,0,0) coordinate (o) -- ++(-\cubescale*\cubex,0,0) coordinate (a) -- ++(0,-\cubescale*\cubey,0) coordinate (b) edge coordinate [pos=1] (g) ++(0,0,-\cubescale*\cubez)  -- ++(\cubescale*\cubex,0,0) coordinate (c) -- cycle
		(o) -- ++(0,0,-\cubescale*\cubez) coordinate (d) -- ++(0,-\cubescale*\cubey,0) coordinate (e) edge (g) -- (c) -- cycle
		(o) -- (a) -- ++(0,0,-\cubescale*\cubez) coordinate (f) edge (g) -- (d) -- cycle;
		\path [every edge/.append style={pic actions, |-|}]
		(b) +(0,-5pt) coordinate (b1) edge ["22 \cubeunits"'] (b1 -| c)
		(b) +(-5pt,0) coordinate (b2) edge ["22 \cubeunits"] (b2 |- a)
		(c) +(3.5pt,-3.5pt) coordinate (c2) edge ["22 \cubeunits"'] ([xshift=3.5pt,yshift=-3.5pt]e)
		;
	},
	/cuboid/.search also={/tikz},
	/cuboid/.cd,
	width/.store in=\cubex,
	height/.store in=\cubey,
	depth/.store in=\cubez,
	units/.store in=\cubeunits,
	scale/.store in=\cubescale,
	width=10,
	height=10,
	depth=10,
	units=cm,
	scale=.1,
}

\usepackage{eso-pic}
\AddToShipoutPicture*{
	\footnotesize\sffamily\raisebox{0.8cm}{\hspace{1.4cm}\fbox{
		\parbox{\textwidth}{
			© 2025 IEEE. Personal use of this material is permitted. Permission from IEEE must be obtained for all other uses, including reprinting/republishing this material for advertising or promotional purposes, collecting new collected works for resale or redistribution to servers or lists, or reuse of any copyrighted component of this work in other works.
		}
	}}
}

\begin{document}

\title{A stabilized Two-Step Formulation of Maxwell's Equations in the time-domain}

\author{\IEEEauthorblockN{Leon Herles\IEEEauthorrefmark{1},
Mario Mally\IEEEauthorrefmark{1,2},
Jörg Ostrowski\IEEEauthorrefmark{3}, and
Sebastian Schöps\IEEEauthorrefmark{1},
Melina Merkel\IEEEauthorrefmark{1}}
\IEEEauthorblockA{\IEEEauthorrefmark{1}Computational Electromagnetics Group, Technische Universität Darmstadt, 64289 Darmstadt, Germany}
\IEEEauthorblockA{\IEEEauthorrefmark{2}Department of Applied Mathematics, Universidade de Santiago de Compostela, 15782 Santiago de Compostela, Spain}
\IEEEauthorblockA{\IEEEauthorrefmark{3}Siemens Digital Industries Software, Switzerland}%
\thanks{%
    Received July 2025; revised August 2025; accepted October 2025. Date of publication October 2025; date of current version October 2025. Corresponding author: L. Herles (e-mail: leon.herles@stud.tu-darmstadt.de).%
}
\thanks{%
    Color versions of one or more figures in this article are available at
https://doi.org/10.1109/TMAG.2025.3619844.%
}
\thanks{%
    Digital Object Identifier 10.1109/TMAG.2025.3619844.%
}
}

\markboth{A stabilized Two-Step Formulation of Maxwell's Equations in the time-domain}%
{A stabilized Two-Step Formulation of Maxwell's Equations in the time-domain}

\IEEEtitleabstractindextext{%
\begin{abstract}
Simulating electromagnetic fields across broad frequency ranges is challenging due to numerical instabilities at low frequencies. This work extends a stabilized two-step formulation of Maxwell’s equations to the time-domain. Using a Galerkin discretization in space, we apply two different time-discretization schemes that are tailored to the first- and second-order in time partial differential equations of the two-step solution procedure used here. To address the low-frequency instability, we incorporate a generalized tree-cotree gauge that removes the singularity of the curl-curl operator, ensuring robustness even in the static limit. Numerical results on academic and application-oriented 3D problems confirm stability, accuracy, and the method’s applicability to nonlinear, temperature-dependent materials.
\end{abstract}

}

\maketitle
\IEEEdisplaynontitleabstractindextext
\IEEEpeerreviewmaketitle

\section{Introduction}
\IEEEPARstart{S}{imulations} across broad frequency ranges are essential for evaluating many electromagnetic devices, particularly in the context of electromagnetic compatibility. This includes static, low-, and high-frequency regimes, where standard formulations are unstable and static or quasistatic approximations may not be sufficiently accurate \cite{Koch_2011aa,Clemens_2021aa,Clemens_2022aa}. 
In the low-frequency range especially, it should also be possible to account efficiently for nonlinear effects, such as field- and/or temperature-dependent material parameters in, e.g., iron or zinc-oxide.

Various stabilization approaches for Darwin and full Maxwell formulations have been proposed, particularly in frequency-domain, e.g. \cite{Dyczij-Edlinger_1999aa,Hiptmair_2008aa,Jochum_2015aa, Ho_2016aa,Eller_2017aa,Badics_2023aa}. However, only few have been translated into the time-domain, e.g. recently the Darwin-based approaches of Zhao \cite{Zhao_2019ac} and Kaimori \cite{Kaimori_2024aa}.
We follow the full Maxwell two-step formulation originally proposed in \cite{Ostrowski_2021aa}. It does not use fractional powers of the frequency and allows a convenient separation of capacitive/resistive and inductive effects in many applications. It was recently low-frequency stabilized using a tree-cotree decomposition in \cite{Herles_2025aa}. In this contribution, we propose a transformation into the time-domain that is stable with respect to large time steps, i.e., when approaching the static limit. Moreover, time discretization can be customized to both steps individually. We use the trapezoidal rule for the first step that has first order derivatives in time, and the Newmark-beta-scheme for the second step that has second order derivatives in time.

Note that the first step decouples from the second step in case of  linear problems, i.e., it can be computed for the full time-interval of interest before the second step is executed.

The paper is structured as follows. In \autoref{sec:two-step} the two-step formulation is recapitulated and discretized in \autoref{sec:discrete}. Its instability as well as stabilization techniques for large time steps are discussed in \autoref{sec:instability}. Then, \autoref{sec:tests} shows the results of a numerical reference implementation.  \autoref{sec:conclusions} closes the paper with some conclusions.

\section{Formulation}\label{sec:two-step}
We follow \cite{Herles_2025aa,Ostrowski_2021aa}, and start from the two-step Maxwell formulation in frequency-domain given by
 \begin{align}
     -\DIV\left(\left(\conductivity + \imNum \omega \permittivity\right)\GRAD\underline{\phifield}\right) &=\imNum\omega\underline{\rhofield}^{s}, \label{eq:strong1_freq}\\
     \CURL\left(\reluctivity\CURL\underline{\Afield}\right) + \left(\imNum \omega \conductivity - \omega^2 \permittivity\right)\underline{\Afield}  &= \underline{\Jfield}^{\mathrm{s}}-\left(\conductivity + \imNum \omega \permittivity\right)\GRAD\underline{\phifield}, \label{eq:strong2_freq}
 \end{align}

on a bounded, open and simply-connected domain $\Omega\subset\mathbb{R}^3$, where $\imNum$ denotes the imaginary unit and $\omega$ the angular frequency. The unknowns are the magnetic vector potential and scalar potential phasors $\underline{\Afield}$ and $\underline{\phifield}$, respectively, that are computed for a given charge density phasor $\underline{\rhofield}^{s}$ and source current density phasor $\underline{\Jfield}^{\mathrm{s}}$. We define the permittivity $\permittivity>0$ and reluctivity $\reluctivity>0$ in $\Omega$ as well as the conductivity $\conductivity>0$ in $\Omega_{\mathrm{C}}\subset\Omega$. In the remaining parts $\overline{\Omega_{\mathrm{A}}}=\overline{\Omega}\setminus\Omega_{\mathrm{C}}$, we assume $\conductivity=0$. 

The corresponding time-domain formulation reads%
\begin{align}
    -\DIV\left(\conductivity \GRAD \phifield + \permittivity \GRAD \dt{\phifield} \right) & = \dt{\rhofield}^{s} \label{eq:strong1_time}\\
    \CURL\left(\reluctivity\CURL\Afield\right) + \conductivity \dt{\Afield} + \permittivity \dtt{\Afield} 
    & = \Jfield^{\mathrm{s}}-\conductivity\GRAD\phifield - \permittivity\GRAD\dt{\phifield} \label{eq:strong2_time}
\end{align}
on an interval $\mathcal{I} = (0,T)$ with given initial values for the fields $\phifield$ and $\Afield$ at $t=0$. In both formulations, \eqref{eq:strong1_freq} and \eqref{eq:strong1_time}, one solves an electroquasistatic (EQS) problem for a given charge density first to obtain the electric scalar potential.
In electroquasistatically dominated applications \eqref{eq:strong2_freq} and \eqref{eq:strong2_time} can be considered as a magnetic correction step, that must only be computed if inductive effects are relevant.
The right-hand-side of \eqref{eq:strong2_time} consists of an electroquasistatic source term and a source current density $\Jfield^{\mathrm{s}}$. 

We employ a combination of homogeneous Dirichlet ($\phifield=0$; $\Afield\times\nfield=\zerofield$) and homogeneous Neumann ($\partial\phifield/\partial\nfield=0$; $\CURL\Afield\times\nfield=\zerofield$) boundary conditions (BCs), where $\nfield$ is the normal vector on the boundary. 

\section{Discretization}\label{sec:discrete}
We employ the method of lines, i.e., we follow the typical procedure of discretizing space to obtain a semi-discrete system and employ time-stepping schemes to compute a corresponding, discrete solution \cite{Kulchytska-Ruchka_2022aa}.

\subsection{Spatial Discretization}
To discretize space, we use a Galerkin approach, for which we omit the corresponding weak formulations since they can be easily derived following the procedure described in \cite{Herles_2025aa}. In the following, we use the discrete subspaces $\scalarDiscreteSpace=\operatorname{span}\{\scalarBaseFun_i\}_{i=1}^{n_{\mathrm{v}}}\subset\scalarWeakSpace$ and $\vectorDiscreteSpace=\operatorname{span}\{\vectorBaseFun_i\}_{i=1}^{n_{\mathrm{w}}}\subset\vectorWeakSpace$, that shall fulfill a discrete De Rham sequence \cite{Monk_2003aa}.
This leads to the matrices
\begin{align}
    \left(\stiffMatGrad{\star}\right)_{ij} &= \int_{\Omega}\star(\GRAD \scalarBaseFun_j)\cdot(\GRAD \scalarBaseFun_i)\,\diffV, \\
    \left(\gradMat{\star}\right)_{ij} &= \int_{\Omega}\star(\GRAD \scalarBaseFun_j)\cdot\vectorBaseFun_i\,\diffV, \label{eq:gradMat} \\
    \left(\massMatCurl{\star}\right)_{ij} &= \int_{\Omega}\star\vectorBaseFun_j\cdot\vectorBaseFun_i\,\diffV,\\
\intertext{where $\star\in\{\conductivity,\permittivity\}$ denotes the employed material property and the curl-curl matrix}
\left(\stiffMatCurl{\reluctivity}\right)_{ij} &= \int_{\Omega}\reluctivity
    (\CURL \vectorBaseFun_j)\cdot(\CURL \vectorBaseFun_i)\,\diffV.
\end{align}
The time dependent source contributions are given by the vectors
\begin{align*}
    \left(\dt{\mathbf{q}}_{\mathrm{s}}\right)_{i} &= \int_{\Omega} \dt{\rhofield}^{s} \scalarBaseFun_i\,\diffV,
    \quad&\text{and}&\quad&
    \left(\mathbf{j}_{\mathrm{s}}\right)_{i} &= \int_{\Omega} \Jfield^{s} \cdot \vectorBaseFun_i\,\diffV.
\end{align*}
Altogether, we end up with the system of ordinary differential equations (ODEs)
\begin{align}
    \stiffMatGrad{\conductivity} \scalarDOFs + \stiffMatGrad{\permittivity} \dt{\scalarDOFs} &= \dt{\mathbf{q}}_{\mathrm{s}} \label{eq:discrete1}\\
    \stiffMatCurl{\reluctivity}\vectorDOFs + \massMatCurl{\conductivity} \dt{\vectorDOFs} + \massMatCurl{\permittivity} \dtt{\vectorDOFs} &= \mathbf{j}(\mathbf{u}, \mathbf{\dt{u}}),\label{eq:discrete2}
\end{align}
in which $\mathbf{j}(\scalarDOFs, \dt{\scalarDOFs}) = \mathbf{j}_{\mathrm{s}} - \gradMat{\conductivity}\scalarDOFs - \gradMat{\permittivity}\dt{\scalarDOFs}$. The degrees of freedom (DOFs) of the scalar and vector potential are time-dependent and given as $\scalarDOFs(t)\in\mathbb{R}^{n_{\mathrm{v}}}$ and $\vectorDOFs(t)\in\mathbb{R}^{n_{\mathrm{w}}}$ for every $t\in\mathcal{I}$. 
\subsection{Temporal Discretization}
Equations \eqref{eq:discrete1} and \eqref{eq:discrete2} can be solved subsequently in time-domain by different approaches. Using methods of the same order in accuracy is recommended to ensure efficiency. For time integration of the first equation \eqref{eq:discrete1}, we apply the trapezoidal rule, which is second order accurate \cite[Chapter II.7]{Hairer_2000aa}. Application of the rule yields
\begin{equation}
\begin{aligned}
    \left(\frac{2}{\deltaT}\stiffMatGrad{\permittivity} + \stiffMatGrad{\conductivity}\right)\scalarDOFs_{n+1} = \left( \frac{2}{\deltaT}\stiffMatGrad{\permittivity} - \stiffMatGrad{\conductivity}\right)\scalarDOFs_{n} \\ + \dt{\mathbf{q}}_{s,n+1} + \dt{\mathbf{q}}_{s,n} \label{eq:trapezoidalEQS}
    \end{aligned}
\end{equation}
with time step size $\deltaT$, where the index $n$ denotes that the solution is approximated at time $n \Delta t$. Let us observe again, that this equation can be readily solved on the whole time interval $\mathcal{I}$.

The ODE system \eqref{eq:discrete2} has second derivatives in time.
We employ the Newmark-beta method \cite{Newmark_1959aa} which is up to second order accurate depending on the choice of parameters $\beta$ and $\gamma$. Application of this methods results in
\begin{align}
    \stiffMatCurl{\mathrm{N}\beta}\vectorDOFs_{n+1} &= \mathbf{f}_{\mathrm{N}\beta}\label{eq:NbetaSys}\\
    \ddot{\vectorDOFs}_{n+1} &= \frac{1}{\deltaT^2 \beta}\left( \vectorDOFs_{n+1} - \vectorDOFs_n - \deltaT \dot{\vectorDOFs}_n\right) - \frac{1-2\beta}{2\beta}\ddot{\vectorDOFs}_n \label{eq:nbAcceleration}\\
    \dot{\vectorDOFs}_{n+1} &= \dot{\vectorDOFs}_n + \left( 1 - \gamma\right) \deltaT \ddot{\vectorDOFs}_n + \gamma \deltaT \ddot{\vectorDOFs}_{n+1} \label{eq:nbVelocity}
\end{align}
with time step size $\deltaT$, update matrix
\begin{equation}
    \stiffMatCurl{N\beta} = \stiffMatCurl{\reluctivity} + \frac{\gamma}{\deltaT \beta} \massMatCurl{\conductivity} \frac{1}{\deltaT^2 \beta} \massMatCurl{\permittivity} \label{eq:updateMatrix}
\end{equation}
and right-hand side
\begin{equation}
\begin{aligned}
     \mathbf{f}_{\mathrm{N}\beta} = \mathbf{j}_{n+1} + \massMatCurl{\permittivity} \left( \frac{1}{\deltaT^2 \beta} \vectorDOFs_n + \frac{1}{\deltaT \beta} \dot{\vectorDOFs}_n + \frac{1-2\beta}{2\beta}\ddot{\vectorDOFs}_n\right) + \\\massMatCurl{\conductivity} \left( \frac{\gamma}{\beta \deltaT} \vectorDOFs_n + \left( \frac{\gamma}{\beta} - 1 \right) \dot{\vectorDOFs}_n + \left( \frac{\gamma}{2 \beta} - 1\right)\deltaT \ddot{\vectorDOFs}_n\right), \label{eq:nbrhs}
\end{aligned}
\end{equation}
where $\mathbf{j}_{n+1}$ is consistently approximated (2nd order) by
\begin{equation}
    \mathbf{j}_{n+1} = \mathbf{j}_{s,n+1} + \frac{1}{2\deltaT}\gradMat{\permittivity}  \left( \scalarDOFs_{n} - \scalarDOFs_{n+2} \right) - \gradMat{\conductivity} \scalarDOFs_{n+1}\;.
\end{equation}
We choose the default parameters $\beta = \frac{1}{4}$ and $\gamma = \frac{1}{2}$ which again corresponds to the trapezoidal rule. Only one linear system \eqref{eq:NbetaSys} has to be solved in every iteration to obtain $\vectorDOFs_{n+1}$ whereas the derivatives can then be computed from \eqref{eq:nbAcceleration} and \eqref{eq:nbVelocity}. Furthermore, this step can be computed after completing \eqref{eq:trapezoidalEQS} on $\mathcal{I}$ or even in parallel (after the first two time steps).

\section{Instability of the Full Maxwell model}
\label{sec:instability}

For large $\deltaT$, i.e., for $\deltaT\rightarrow\infty$, one can see that most terms vanish in \eqref{eq:trapezoidalEQS}. The remaining system matrix is usually singular as $\sigma$ vanishes in the nonconducting parts of the domain $\Omega_\mathrm{A}$, thus leading to stability issues. While this first step is not the focus of this contribution, effective modifications are proposed in \cite{Balian_2023aa}. This instability for large time steps $\deltaT$ matches the low frequency breakdown occurring in the frequency-domain when approaching the static limit, i.e. for vanishing angular frequency $\omega \rightarrow 0$. This is due to the equivalence of $\imNum \omega$ in the frequency-domain with ${1}/{\deltaT}$ in the time-domain \cite{Ostrowski_2012aa}.

The second step \eqref{eq:updateMatrix} is haunted by a similar instability for large $\deltaT$. One can see that only the singular $\stiffMatCurl{\reluctivity}$ remains in \eqref{eq:updateMatrix} due to the division by $\deltaT$. Its kernel consists of discrete gradient fields, and additional issues arise because of large scaling differences due to $ \frac{\conductivity}{\deltaT}\gg\frac{\permittivity}{\deltaT^2}$.

Generalizing the ideas of \cite{Herles_2025aa} from frequency to time-domain, we split the constraint
\begin{equation}
    \DIV\left(\conductivity\dt{\Afield} + \permittivity \dtt{\Afield} \right) = 0, \label{eq:strongGenCoul}
\end{equation}
into the corresponding subdomains $\Omega_\mathrm{C}$ and $\Omega_\mathrm{A}$ and then discretize it using the same Galerkin approach as before and integrate the element-wise equations in time such that we obtain
\begin{equation}
    \underbrace{\begin{bmatrix}
        \divMat{\conductivity}^{(\mathrm{CC})} & \divMat{\conductivity}^{(\mathrm{CA})} \\
        \lambda\divMat{\permittivity}^{(\mathrm{AC})} & \lambda\divMat{\permittivity}^{(\mathrm{AA})} \\
    \end{bmatrix}}_{\textstyle = \tilde{\divMat{}}}
    \vectorDOFs 
    + \underbrace{\begin{bmatrix}
        \divMat{\permittivity}^{(\mathrm{CC})} & \divMat{\permittivity}^{(\mathrm{CA})} \\
        \mathbf{0} & \mathbf{0}  \\
    \end{bmatrix}}_{\textstyle =\hat{\divMat{}}} 
    \dot{\vectorDOFs}
    = 
    \begin{bmatrix}
        \mathbf{0} \\
        \mathbf{0}
    \end{bmatrix}
    ,
    \label{eq:stableDiv}
\end{equation}
with superscripts $\mathrm{C}, \mathrm{A}$ denoting the support of the basis functions of the respective parts of the domain $\Omega_\mathrm{C}, \Omega_\mathrm{A}$, a scaling factor $\lambda$ and the weighted divergence matrix
\begin{equation}
    \left(\divMat{\star}\right)_{ij} = \int_{\Omega}\star \DIV\left(\vectorBaseFun_j\right)  \scalarBaseFun_i\,\diffV \label{eq:divMat}
\end{equation}
for $\star\in\{\conductivity,\permittivity\}$. The scaling factor can be used to improve matrix conditioning, let us assume $\lambda = 1$ for now. Note that the reformulation \eqref{eq:stableDiv} remains stable even if $\dot{\vectorDOFs}\to0$ (or $\deltaT\to\infty$) since $\tilde{\divMat{}}$ is invertible.
This constraint can now be conveniently included in \eqref{eq:NbetaSys} for stabilization.

\subsection{Tree-Cotree Decomposition}
To eliminate the discrete kernel of the stiffness matrix $\stiffMatCurl{\reluctivity}$, we employ a discrete tree-cotree gauge \cite{Albanese_1988aa,Munteanu_2002aa,Manges_1995aa,Eller_2017aa}, which relies solely on topological mesh information. This approach interprets the mesh as a graph and decomposes it into a spanning tree and a remaining cotree. This decomposition allows the identification of determined and underdetermined DOFs, where the underdetermined DOFs, denoted by $\vectorDOFs^{(\mathrm{T})}$, correspond to the edges of the spanning tree. The remaining DOFs $\vectorDOFs^{(\mathrm{R})}$, associated with the cotree edges, are then determined from the magnetostatic problem
\begin{equation}
    \begin{bmatrix}
	    \stiffMatCurl{\reluctivity}^{(\mathrm{RR})} & \stiffMatCurl{\reluctivity}^{(\mathrm{RT})}\\
	    \stiffMatCurl{\reluctivity}^{(\mathrm{TR})} & \stiffMatCurl{\reluctivity}^{(\mathrm{TT})}
	\end{bmatrix}
	\begin{bmatrix}
	    \vectorDOFs^{(\mathrm{R})} \\
	    \vectorDOFs^{(\mathrm{T})}
	\end{bmatrix}
	=
	\begin{bmatrix}
	\mathbf{j}^{(\mathrm{R})}(\mathbf{u})\\
	\mathbf{j}^{(\mathrm{T})}(\mathbf{u})
	\end{bmatrix}. \label{eq:splitMagneto}
\end{equation}
In the simplest form of tree-cotree gauging proposed in \cite{Albanese_1988aa} the tree DOFs are set to zero, i.e., $\vectorDOFs^{(\mathrm{T})}=\mathbf{0}$. Alternatively, the underdetermined can be constrained by prescribing a relation to the remaining DOFs \cite{Rapetti_2022aa,Munteanu_2002aa}, e.g.,
\begin{equation}
    \vectorDOFs^{(\mathrm{T})}=\stiffMatCurl{\reluctivity}^{(\mathrm{TR})}\left(\stiffMatCurl{\reluctivity}^{(\mathrm{RR})}\right)^{-1}\vectorDOFs^{(\mathrm{R})}.\label{eq:discrOrtho}
\end{equation}
Then, $\vectorDOFs^{(\mathrm{R})}$ is computed by reformulating \eqref{eq:splitMagneto}. Note that $\stiffMatCurl{\reluctivity}^{(\mathrm{RR})}$ is an invertible matrix and that for any given $\vectorDOFs^{(\mathrm{T})}$, the $\vectorDOFs^{(\mathrm{R})}$ obtained from the first row of \eqref{eq:splitMagneto}, automatically  satisfies the second row. 
Note that the tree-cotree decomposition is performed globally on the domain $\Omega$ and the resulting gauging approach is applied consistently in both the conducting and nonconducting subdomains, $\Omega_{\mathrm{C}}$ and $\Omega_{\mathrm{A}}$, respectively.

\subsection{Stabilization}
Following the approach outlined in \cite{Herles_2025aa}, we can reorder the DOFs in \eqref{eq:NbetaSys} similarly as in \eqref{eq:splitMagneto} to obtain
\begin{align}
    \stiffMatCurl{N\beta} 
    = 
    &
    \begin{bmatrix}
	    \stiffMatCurl{\reluctivity}^{(\mathrm{RR})} & \stiffMatCurl{\reluctivity}^{(\mathrm{RT})}\\
	    \stiffMatCurl{\reluctivity}^{(\mathrm{TR})} & \stiffMatCurl{\reluctivity}^{(\mathrm{TT})}
	\end{bmatrix} 
    +
    \frac{\gamma}{\deltaT \beta}
    \begin{bmatrix}
	    \massMatCurl{\conductivity}^{(\mathrm{RR})} & \massMatCurl{\conductivity}^{(\mathrm{RT})} \\
	    \massMatCurl{\conductivity}^{(\mathrm{TR})} & \massMatCurl{\conductivity}^{(\mathrm{TT})}
	\end{bmatrix} 
    \nonumber
    \\ 
    &
    +
    \frac{1}{\deltaT^2 \beta}
    \begin{bmatrix}
	    \massMatCurl{\permittivity}^{(\mathrm{RR})} & \massMatCurl{\permittivity}^{(\mathrm{RT})}\\
	    \massMatCurl{\permittivity}^{(\mathrm{TR})} & \massMatCurl{\permittivity}^{(\mathrm{TT})}
	\end{bmatrix}. \label{eq:updateMatrixSplit}
\end{align}
We now replace the lower blocks of equations in \eqref{eq:updateMatrixSplit} 
which become linearly dependent if $\deltaT \rightarrow \infty$ with the stabilized version of the implicit gauge condition \eqref{eq:stableDiv}. Finally, we also reorder \eqref{eq:stableDiv} to arrive at
\begin{equation}
    \begin{bmatrix}
        \tilde{\divMat{}}^{(\mathrm{R})} &
        \tilde{\divMat{}}^{(\mathrm{T})} 
    \end{bmatrix}
    \begin{bmatrix}
	    \vectorDOFs^{(\mathrm{R})} \\
	    \vectorDOFs^{(\mathrm{T})}
	\end{bmatrix} 
    + 
    \begin{bmatrix}
        \hat{\divMat{}}^{(\mathrm{R})} &
        \hat{\divMat{}}^{(\mathrm{T})} 
    \end{bmatrix}     \begin{bmatrix}
	    \dot{\vectorDOFs}^{(\mathrm{R})} \\
	    \dot{\vectorDOFs}^{(\mathrm{T})}
	\end{bmatrix} 
    =
    0,
    \label{eq:stableDivTC}
\end{equation}
and then combine it with the first line of \eqref{eq:updateMatrixSplit} to obtain
\begin{align}
    \stiffMatCurl{N\beta}
    =&
    \begin{bmatrix}
	    \stiffMatCurl{\reluctivity}^{(\mathrm{RR})} & \stiffMatCurl{\reluctivity}^{(\mathrm{RT})}\\
	    \tilde{\divMat{}}^{(\mathrm{TR})} & \tilde{\divMat{}}^{(\mathrm{TT})}
	\end{bmatrix}
    +
    \frac{\gamma}{\deltaT \beta} \begin{bmatrix}
	    \massMatCurl{\conductivity}^{(\mathrm{RR})} & \massMatCurl{\conductivity}^{(\mathrm{RT})} \\
	    \hat{\divMat{}}^{(\mathrm{TR})} & \hat{\divMat{}}^{(\mathrm{TT})}
	\end{bmatrix}
    \nonumber
    \\ 
    &+
    \frac{1}{\deltaT^2 \beta}
    \begin{bmatrix}
	    \massMatCurl{\permittivity}^{(\mathrm{RR})} & \massMatCurl{\permittivity}^{(\mathrm{RT})}\\
	    \mathbf{0} & \mathbf{0}
	\end{bmatrix}, 
    \label{eq:updateMatrixStable}
\end{align}
which does not contain the kernel of $\stiffMatCurl{\reluctivity}$ anymore. Computation of the derivatives is still possible sequentially using \eqref{eq:nbAcceleration} and \eqref{eq:nbVelocity}. The right-hand side vector \eqref{eq:nbrhs} also has to be sorted and decomposed accordingly.

\section{Numerical Tests}\label{sec:tests}
We consider two 3D problems for the evaluation of the performance of the proposed stabilization method. The first problem is an academic example to investigate the numerical behavior of the method, while the second example deals with a more application-oriented configuration. We employ \texttt{GeoPDEs} \cite{Vazquez_2016aa} for our numerical experiments. The scaling factor for the divergence matrix in \eqref{eq:stableDiv} is chosen as
\begin{equation}
    \lambda = \frac{\max_{\boldsymbol{x} \in \Omega} \conductivity (\boldsymbol{x}) + \conductivity_\mathrm{art}}{\max_{\boldsymbol{x} \in \Omega } \varepsilon(\boldsymbol{x})},
\end{equation}
to improve conditioning of the resulting linear systems where we set $\conductivity_\mathrm{art} = \SI{1e-6}{}$ to ensure $\lambda \neq 0$.

\subsection{Academic Example}
The first test problem consists of three conducting bars in a dielectric box. Its configuration is shown in \autoref{fig:academicTestEx}. It is inspired by the test problem of \cite{Ostrowski_2021aa}. As for boundary conditions, we employ  a homogeneous Dirichlet boundary $\Afield \times \nfield = \zerofield$ on the whole boundary $\partial\Omega$ for the magnetic vector potential. The boundary conditions for the electric scalar potential are chosen in such a way that the resulting solution of the electric flux density $\Dfield_{\mathrm{e}} = -\permittivity \GRAD \phifield$ of the first step, i.e., the electroquasistatic problem, is constant and known in closed form. The second step then computes the inductive correction step corresponding to the skin effect in this case. The excitation is realized via a sinusoidal voltage excitation 
\begin{equation}
    V_1 =  \hat{u} \sin{2\pi f t},
    \label{eq:Excitation}
\end{equation}
with $\hat{u} = \SI{1}{\volt}$ and $ f= \SI{150}{Hz}$. For both steps, we employ zero initial conditions. The problem is discretized in space using 27 patches and 6084 third order basis functions. We discretize it in time using 20 points per period and simulate two periods of the excitation for a total time of $\mathcal{I} = (0, \SI{13.3}{ms})$ resulting in a time step size of $\deltaT = \SI{0.333}{ms}$.

\begin{figure}
    \centering
    \begin{subfigure}[B]{0.32\linewidth}
        \centering
        \includegraphics[width=0.9\linewidth]{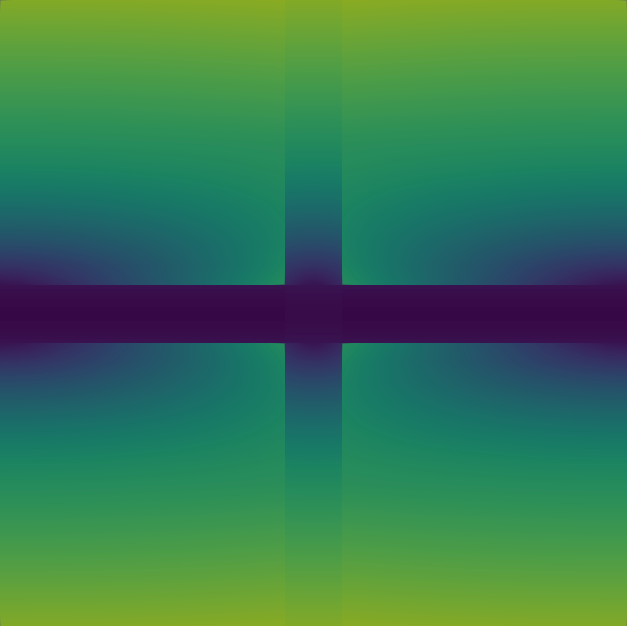}
        \caption{$t_1=\SI{0.17}{ms}$.}
        \label{fig:Dfield_1}
    \end{subfigure}
    \begin{subfigure}[B]{0.32\linewidth}
        \centering
        \includegraphics[width=0.9\linewidth]{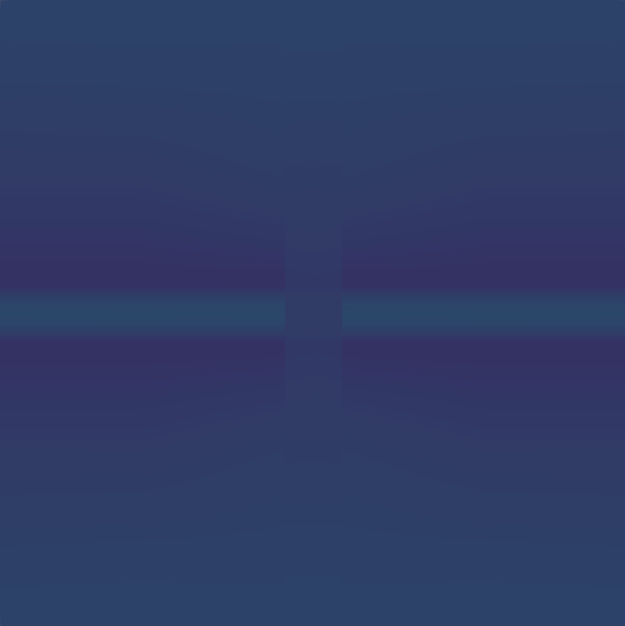}
        \caption{$t_2=\SI{0.34}{ms}$.}
        \label{fig:Dfield_100}
    \end{subfigure}
    \begin{subfigure}[B]{0.32\linewidth}
        \centering
        \includegraphics[width=0.9\linewidth]{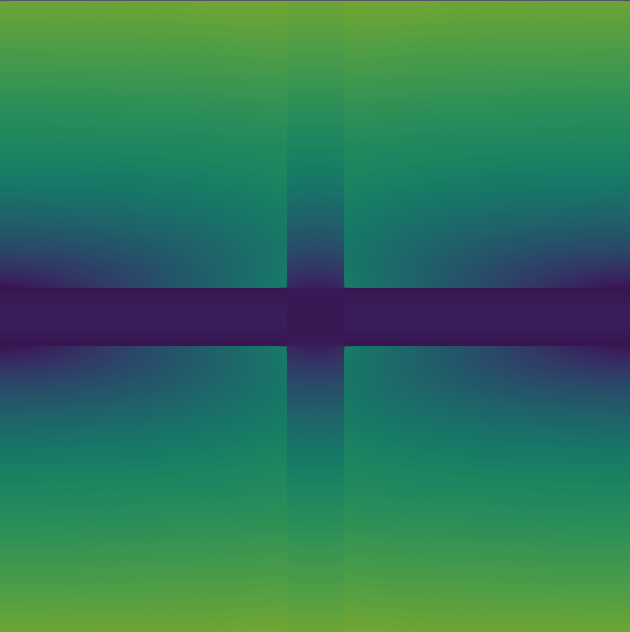}
        \caption{$t_3=\SI{11.7}{ms}$.}
        \label{fig:Dfield_1000}
    \end{subfigure}
     \begin{subfigure}[B]{\linewidth}
        \centering
        \begin{tikzpicture}
        \pgfplotsset{every tick label/.append style={font=\small}}
            \begin{axis}[%
  hide axis,
  scale only axis,
  colorbar/width=2mm,
  width = 50mm,
  anchor=south,
  point meta min=0.0,
  point meta max=1.7e-10,
  colormap/viridis,                     %
  colorbar horizontal,                  %
  colorbar sampled,                     %
  colorbar style={
    separate axis lines,
    samples=256,                        %
    xlabel=\small{Electric flux density $\big(\SI{}{\coulomb/\meter^{2}}\big)$},
    xtick={0.0, 0.85e-10, 1.7e-10},
    scaled x ticks=false
  },
]
  \addplot [draw=none] coordinates {(0,0)};
\end{axis}
        \end{tikzpicture}
    \end{subfigure}
  \caption{Electric flux density $\norm{\Dfield}{2}$ for different times.}
  \label{fig:academicFields}
\end{figure}

\autoref{fig:academicFields} shows $\norm{\Dfield}{2}$ for different time steps. In \autoref{fig:Dfield_1} one can see that the field for $t_1 = \SI{1.7}{ms}$ is being pulled out of the conductor due to the skin effect. It also corresponds well with the frequency based solution computed in \cite{Herles_2025aa}. This behavior is expected as $t_1$ corresponds to the peak value of the excitation. \autoref{fig:Dfield_100} and \autoref{fig:Dfield_1000} show the resulting flux density for $t_2 = \SI{3.4}{ms}$ and $t_3 = \SI{11.7}{ms}$ respectively. Time $t_2$ corresponds to the zero crossing of the excitation and thus contains small fields, while $t_3$ corresponds again to the (negative) peak value of the excitation and thus matches the flux density of $t_1$ in magnitude. Simulation of this problem using the original (unstabilized) method leads to diverging behavior. This is again in agreement with the behavior observed in \cite{Herles_2025aa}.

\def\LEN{5}
\begin{figure}
    \centering
      \subfloat[Material distribution.\label{fig:Cube_mat}]{%
       \scalebox{0.53}{
                    \begin{tikzpicture}
                        \tikzset{every node/.append style = {font=\large}}
                		\draw[black, very thick] (0,0) rectangle (\LEN,\LEN);

                		\draw (0,-1) -- (\LEN,-1);
                		\node[xshift={3mm}] at (\LEN*2/15,-0.7) {10 cm};
                		\node[xshift={-1mm}] at (\LEN*6/11,-0.7) {2 cm};
                		\node[xshift={-5mm}] at (\LEN*14/15,-0.7) {10 cm};
                		\draw (0,-1.1) -- (0,-0.9);
                		\draw (\LEN*10/22,-1.1) -- (\LEN*10/22,-0.9);
                		\draw (\LEN*12/22,-1.1) -- (\LEN*12/22,-0.9);
                		\draw (\LEN,-1.1) -- (\LEN,-0.9);

                		\draw[fill=TUDa-1b,fill opacity=0.2] (0,0) rectangle (\LEN*10/22,\LEN*10/22);
                        \draw[fill=TUDa-1b,fill opacity=0.2] (\LEN*12/22,0) rectangle (\LEN,\LEN*10/22);
                        \draw[fill=TUDa-1b,fill opacity=0.2] (0,\LEN*12/22) rectangle (\LEN*10/22,\LEN);
                        \draw[fill=TUDa-1b,fill opacity=0.2] (\LEN*12/22,\LEN*12/22) rectangle (\LEN,\LEN);

                        \draw[fill=TUDa-9d, fill opacity=0.2] (\LEN*10/22,0) rectangle (\LEN*12/22,\LEN*10/22);
                        \draw[fill=TUDa-9d, fill opacity=0.2] (\LEN*10/22,\LEN*12/22) rectangle (\LEN*12/22,\LEN);

                		\draw[fill=TUDa-4d, fill opacity=0.2] (0,\LEN*10/22) rectangle (\LEN*10/22,\LEN*12/22);
                		\draw[fill=TUDa-3b, fill opacity=0.5] (\LEN*10/22,\LEN*10/22) rectangle (\LEN*12/22,\LEN*12/22);
                		\draw[fill=TUDa-4d, fill opacity=0.2] (\LEN*12/22,\LEN*10/22) rectangle (\LEN,\LEN*12/22);

                        \node[fill=TUDa-1b, rectangle, inner sep= 0pt, minimum width=5mm, minimum height=5mm, fill opacity=0.2](l1) at (\LEN+1,\LEN-0.5){};
                        \node[fill=TUDa-9d, rectangle, inner sep= 0pt, minimum width=5mm, minimum height=5mm, fill opacity=0.2,anchor=north](l2) at ($(l1.south)-(0,0.8)$){};
                        \node[fill=TUDa-4d, rectangle, inner sep= 0pt, minimum width=5mm, minimum height=5mm, fill opacity=0.2,anchor=north](l3) at ($(l2.south)-(0,0.8)$){};
                        \node[fill=TUDa-3b, rectangle, inner sep= 0pt, minimum width=5mm, minimum height=5mm, fill opacity=0.5,anchor=north](l4) at ($(l3.south)-(0,0.8)$){};

                        \node[anchor=west] at (l1.north east){$\permittivity_\text{r}\,=5$};
                        \node[anchor=west] at (l1.south east){$\conductivity_\text{r}=0$};
                        \node[anchor=west] at (l2.north east){$\permittivity_\text{r}\,=1$};
                        \node[anchor=west] at (l2.south east){$\conductivity_\text{r}=0$};
                        \node[anchor=west] at (l3.north east){$\permittivity_\text{r}\,=5$};
                        \node[anchor=west] at (l3.south east){$\conductivity_\text{r}=5$};
                        \node[anchor=west] at (l4.north east){$\permittivity_\text{r}\,=1$};
                        \node[anchor=west] at (l4.south east){$\conductivity_\text{r}=1$};

                        \draw[black, very thick] (0,\LEN*10/22) rectangle (\LEN,\LEN*12/22);

                		\draw[-{Stealth[black]}] (0,-1.5)   -- node[pos=1.35, xshift={-3mm}, anchor=west]{$\boldsymbol{e}_z$} (1,-1.5);
                	\end{tikzpicture}
                }
       }
    \hfill
  \subfloat[3D-view.\label{fig:Cube3D}]{%
       \scalebox{0.53}{%
                    \begin{tikzpicture}
                        \tikzset{every node/.append style = {font=\large}}

                        \pic [fill=TUDa-1b!25, text=TUDa-1b, draw=black,opacity=0.5] at (-4.4,0) {cuboid={width=0.1, height=44, depth=44}};
                		\pic [fill=TUDa-4d!20, text=TUDa-1b, draw=black] at (+0.866,-2+0.866) {cuboid={width=20, height=4, depth=4}};
                		\pic [fill=TUDa-4d!20, text=TUDa-1b, draw=black] at (+0.866-2-0.4,-2+0.866) {cuboid={width=20, height=4, depth=4}};
                		\pic [fill=none, text=black, draw=black] at (0,0) {annotated cuboid={width=44, height=44, depth=44}};
                		\pic [fill=TUDa-3b!50, text=TUDa-1b, draw=black] at (+0.866-2,-2+0.866) {cuboid={width=4, height=4, depth=4}};

                		\pic [fill=TUDa-9d!25, text=TUDa-1b, draw=black,opacity=0.5] at (0,0) {cuboid={width=0.1, height=44, depth=44}};

                		\pic [fill=gray, text=TUDa-1b, draw=black, opacity=0.1] at (-2.4,0) {cuboid={width=0, height=44, depth=44}};
                		\pic [fill=gray, text=TUDa-1b, draw=black, opacity=0.1] at (-2,0) {cuboid={width=0, height=44, depth=44}};

                		\draw[-{Stealth[black]}] (0-\LEN,-0.25-\LEN)   -- node[pos=1.35, xshift={-3mm}, anchor=west]{$\boldsymbol{e}_z$} (1-\LEN,-0.25-\LEN);

                		\node[xshift={14mm},yshift={29mm},rotate=45] at (-\LEN-0.3,-\LEN*0.5+0.866) {$\phifield = 0$};
                		\node[xshift={-2mm},rotate=45] at (+0.866 + 0.2,-\LEN*0.5+0.866-0.3) {$\phifield = V_1$};
                		\node[xshift={3mm}] at (-\LEN*0.5+0.866,+0.866) {$\frac{\partial \phifield}{\partial \nfield} = 0$};
                	\end{tikzpicture}
                }
       }
  \caption{Setup of the academic test example; adapted from~\cite{Herles_2025aa}.}
  \label{fig:academicTestEx}
\end{figure}
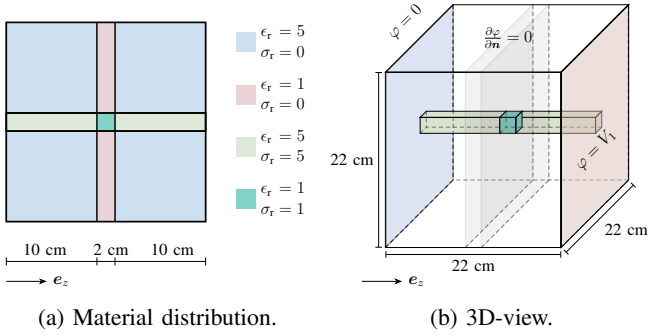

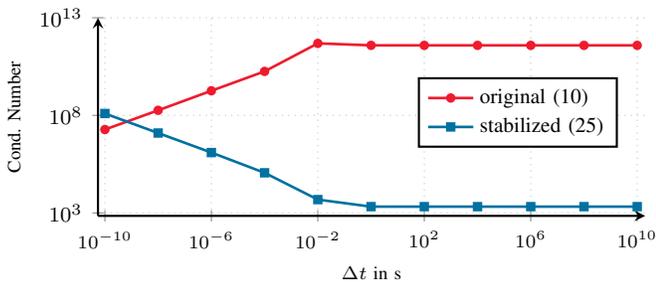
\begin{figure}
    \centering
 \pgfplotstableread[col sep=comma]{data/condOverDt.csv}\cond
                \begin{tikzpicture}
                	\begin{loglogaxis}[tudalineplot, width=\linewidth, height=12em, xlabel={\scriptsize$\deltaT$ in s}, ylabel={\scriptsize Cond. Number}, xticklabel style={font=\scriptsize}, xmin=5.5e-11,xmax=2e10, ymin=5e2, ymax=1e23, xtick distance=10^(4), yticklabels={$10^3$, $10^8$, $10^{13}$, $10^{18}$, $10^{23}$}, xlabel shift={-0mm}, ylabel shift={-0mm},legend style={at={(0.95,0.7)},anchor=north east}]
                    
                 \addplot+[mark size=1.5mm,mark=*,mark options={fill=TUDa-9b,scale=0.3}, TUDa-9b, line width=1pt] table[x index = 0, y index = 1] {\cond};
                 \addplot+[mark size=1.5mm,mark=square*,mark options={fill=TUDa-2c,scale=0.3}, TUDa-2c, line width=1pt] table[x index = 0, y index = 2] {\cond};
                \addlegendentry{original \eqref{eq:discrete2}};
                \addlegendentry{stabilized \eqref{eq:updateMatrixStable}};
                \end{loglogaxis}
                \end{tikzpicture}

                \caption{Condition number of different system matrices of academic test example from \autoref{fig:academicTestEx} over time step size.}
    \label{fig:condNums}
\end{figure}

\autoref{fig:condNums} shows that for large time steps, the original system matrix in \eqref{eq:discrete2} is ill-conditioned or even singular, while the condition number of the stabilized method \eqref{eq:updateMatrixStable} is relatively small for $\deltaT \rightarrow \infty$. 
In the static limit $\frac{1}{\deltaT} = 0 $, we obtain a condition number of approximately $\SI{5.25e5}{}$, while the original update matrix is singular. The condition number of the stabilized system increases slightly in comparison to the one of the original for very small time steps. This effect is well known, see for example \cite{Hiptmair_2000aa} and \cite{Munteanu_2002aa}.

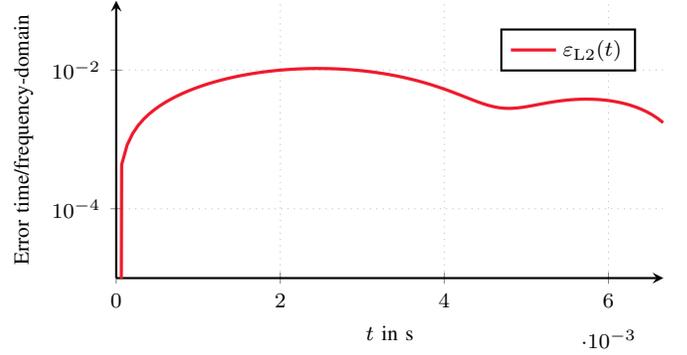
\begin{figure}
\centering
\pgfplotstableread[col sep=comma]{data/cube100.csv}\err
	\begin{tikzpicture}
		\begin{semilogyaxis}[tudalineplot, width=\linewidth, height=15em,
			xlabel={\footnotesize$t$ in s}, ylabel={\footnotesize Error time/frequency-domain}, xticklabel style={font=\footnotesize}, xmin=0,xmax=0.00667, ymin=1e-7,ymax=1e-3, xtick distance=0.002, xlabel shift={-0mm}, ylabel shift={-0mm}, 
            legend style={at={(0.95,0.9)},
            anchor=north east}]	
			\addplot+[mark=none,TUDa-9b, line width=1.2pt] table[x index = 0, y index = 1] {\err};
			\addlegendentry{$\LTwoErr(t)$};		
		\end{semilogyaxis}
\end{tikzpicture}
 \caption{Error \eqref{eq:error} of time domain w.r.t. frequency-domain.}
    \label{fig:FDTDcompare}
\end{figure}

As another method of verification, we compare the solution of the electric field $\Efield_h$ computed in the time-domain with a solution $\underbar{\Efield}_h$ computed in the frequency-domain which is then transformed back to the time-domain using the inverse Fourier transform $\mathcal{F}^{-1}$. The computation in the frequency-domain is carried out using the approach proposed in \cite{Herles_2025aa}. We simulate again the academic cube example from \autoref{fig:academicTestEx} but this time for only one period with 100 time steps to obtain an adequate time-domain precision. The initial value is chosen according to the frequency domain solution. \autoref{fig:FDTDcompare} shows the relative L2-error 
\begin{equation}
    \LTwoErr^2(t) = \frac{\int_\Omega \left( \Efield_h(t) -  \mathcal{F}^{-1}\left\{\underbar{\Efield}_h(\omega)\right\} \right)^2 \diffV}{\max\limits_{t\in\mathcal{I}} \int_\Omega \mathcal{F}^{-1}\left\{\underbar{\Efield}_h(\omega)\right\}^2\diffV}\;.
    \label{eq:error}
\end{equation}
The error is, in the worst case, around 0.01\%. We conclude that both solutions are in good agreement.

 \subsection{Simulation of Setup with Planar Coil} \label{sec:numPlanarCoil}
\begin{figure}
    \centering
    \begin{tikzpicture}
        \node[inner sep=0pt] (test) at (0,0) {\includegraphics[width=\linewidth]{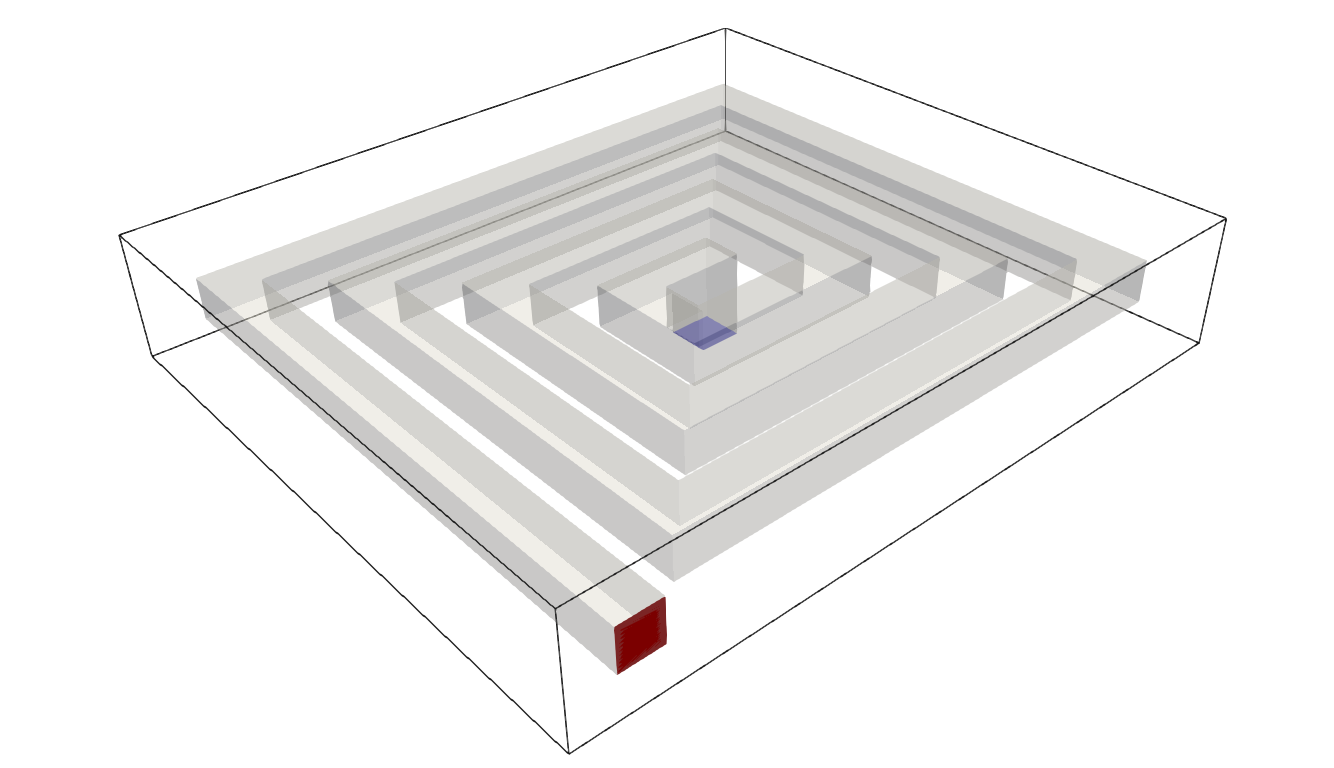}};
        \node[inner sep=0pt] (test) at (-0.2,0) {$\Gamma_\mathrm{G}$};
        \node[inner sep=0pt] (test) at (0.4,-1.5) {$\Gamma_\mathrm{E}$};
    \end{tikzpicture}
    \caption{Problem setup of the planar coil. Taken from \cite{Herles_2025aa}.} 
    \label{fig:planarCoil}
\end{figure}

 Next, we consider a more application-motivated example which is based on a planar coil. The problem setup is shown in \autoref{fig:planarCoil}. The coil is made from copper and has a cross-section of $\SI{3}{mm} \times \SI{3}{mm}$ and a conductivity of $\conductivity = \SI{6e7}{\frac{S}{m}}$. It consists of three {turns} and is surrounded by air. The airbox has dimensions $\SI{45}{mm} \times \SI{51}{mm} \times \SI{9}{mm}$. The excitation is again realized via a sinusoidal voltage \eqref{eq:Excitation}, resulting in boundary conditions $\phifield=\SI{0}{V}$ on $\Gamma_{\mathrm{G}}$ and $\phifield=V_1$ on $\Gamma_\mathrm{E}$. On the remaining boundary $\partial\Omega\setminus \left(\Gamma_\mathrm{G}\cup\Gamma_\mathrm{E}\right)$, we set $\GRAD\phifield\cdot\nfield=0$. For the second step, perfect magnetic boundary conditions are assumed, i.e. $\CURL\Afield\times\nfield=\zerofield$. For both steps, we employ zero initial conditions. The problem is discretized using 144 patches and 12,822 second order basis functions. It is simulated for one full period of the excitation using 20 time steps, resulting in a total simulation time interval of $\mathcal{I} = (0, \SI{6.65}{ms})$ with a time step size of $\deltaT = \SI{0.333}{ms}$. 

\begin{figure}
     \begin{subfigure}[B]{0.48\linewidth}
        \centering
        \includegraphics[width=0.9\linewidth,clip]{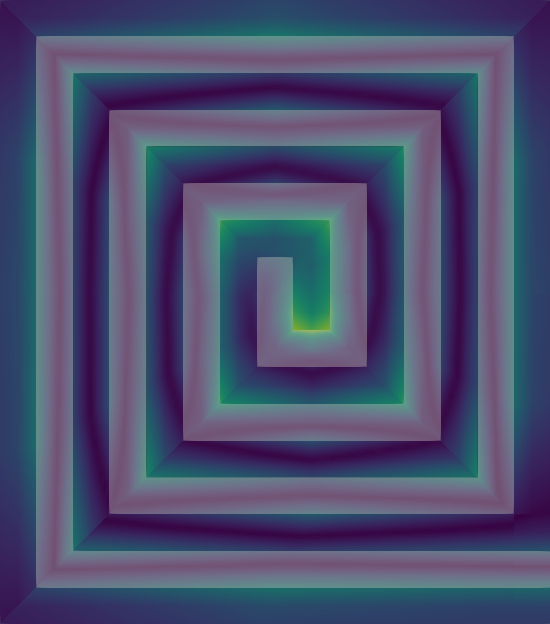}\\
        \hspace{1mm}
    \begin{tikzpicture}[every node/.style={scale=0.8}]
        \pgfplotsset{every tick label/.append style={font=\small}}
            \begin{axis}[%
  hide axis,
  scale only axis,
  colorbar/width=2mm,
  width = 33mm,
  anchor=south,
  point meta min=0.0,
  point meta max=0.42,
  colormap/viridis,                     %
  colorbar horizontal,                  %
  colorbar sampled,                     %
  colorbar style={
    separate axis lines,
    samples=256,                        %
    xlabel=\small{Magnetic flux density $\big(\SI{}{\tesla}\big)$},
    xtick={0.0,0.21, 0.42},
    scaled x ticks=false
  },
]
  \addplot [draw=none] coordinates {(0,0)};
\end{axis}
        \end{tikzpicture}
        \caption{$\norm{\Bfield}{2}$, $t = \SI{1.7}{ms}$}
        \label{fig:Bfield_pC}
    \end{subfigure} 
    \begin{subfigure}[B]{0.48\linewidth}
        \centering
        \includegraphics[width=0.9\linewidth,clip]{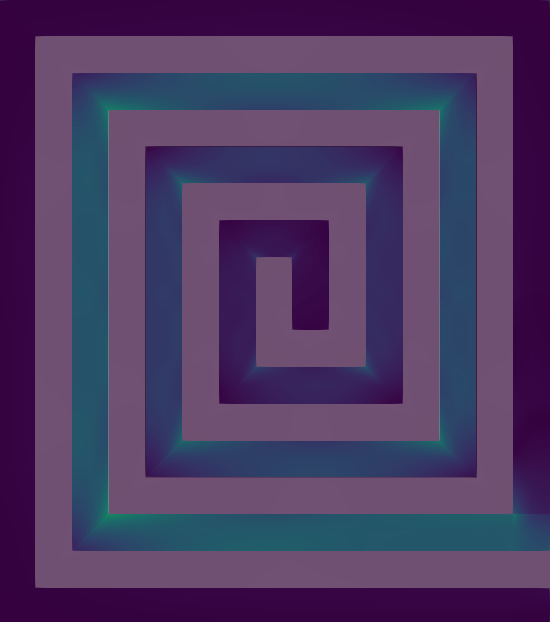}\\
        \hspace{1mm}
        \begin{tikzpicture}[every node/.style={scale=0.8}]
        \pgfplotsset{every tick label/.append style={font=\small}}
            \begin{axis}[%
  hide axis,
  scale only axis,
  colorbar/width=2mm,
  width = 33mm,
  anchor=south,
  point meta min=0.0,
  point meta max=400,
  colormap/viridis,                     %
  colorbar horizontal,                  %
  colorbar sampled,                     %
  colorbar style={
    separate axis lines,
    samples=256,                        %
    xlabel={\small Electric field strength $\big(\SI{}{\volt/\meter}\big)$},
        xtick={0.0,200, 400},
    scaled x ticks=false
  },
]
  \addplot [draw=none] coordinates {(0,0)};
\end{axis}
        \end{tikzpicture}
        \caption{$\norm{\Efield}{2}$, $t = \SI{1.7}{ms}$}
        \label{fig:Efield_pC}
    \end{subfigure}
    \caption{Top view of the field solutions at $t = \SI{1.7}{ms}$ of the planar coil, sliced in the middle.}
    \label{fig:planarCoilFields}
\end{figure}
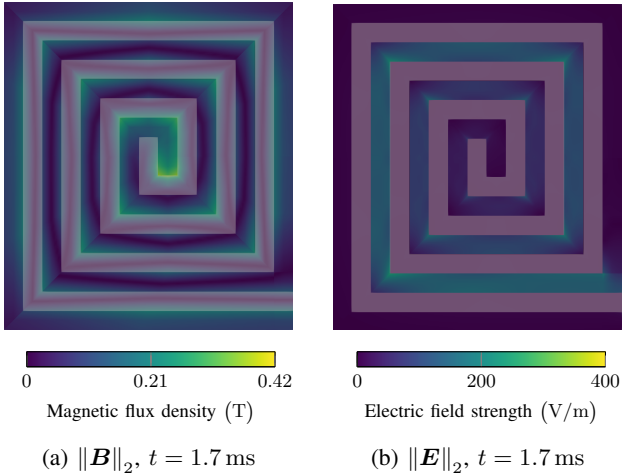

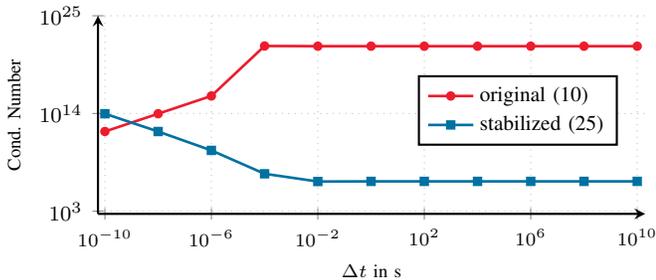
\begin{figure}
 \centering
 \pgfplotstableread[col sep=comma]{data/condOverDtSpuleTest.csv}\condSpule
                \begin{tikzpicture}
                \begin{loglogaxis}[tudalineplot, width=\linewidth, height=12em, xlabel={\scriptsize$\deltaT$ in s}, ylabel={\scriptsize Cond. Number}, xticklabel style={font=\scriptsize}, xmin=5.5e-11,xmax=2e10, ymin=5e2, ymax=1e25,legend style={at={(0.95,0.7)},anchor=north east}]

                 \addplot+[mark size=1.5mm,mark=*,mark options={fill=TUDa-9b,scale=0.3}, TUDa-9b, line width=1pt] table[x index = 0, y index = 1] {\condSpule};
                 \addplot+[mark size=1.5mm,mark=square*,mark options={fill=TUDa-2c,scale=0.3}, TUDa-2c, line width=1pt] table[x index = 0, y index = 2] {\condSpule};
                \addlegendentry{original \eqref{eq:discrete2}};
                \addlegendentry{stabilized \eqref{eq:updateMatrixStable}};
                \end{loglogaxis}
                \end{tikzpicture}
                \caption{Condition number of different system matrices of planar coil over time step size.}
    \label{fig:condNumsSpule}
\end{figure}

In \autoref{fig:planarCoilFields} the magnetic flux density and electric field strength are shown. One can see that the two-step formulation can successfully capture both the capacitive coupling between the turns of the planar coil and the inductive effects of the surrounding magnetic flux density. The coil is highlighted in gray. Note that here again, the nonstable original formulation would lead to nonphysical fields due to badly conditioned update matrices.

\autoref{fig:condNumsSpule} shows the condition numbers of the different methods over the time step size. 
We can again observe that the original system becomes singular for larger time steps. The condition number of the stabilized system is much lower in comparison to the original system and does not deteriorate for $\deltaT \rightarrow \infty$. In the static limit $\frac{1}{\deltaT} = 0$, we obtain a condition number of approximately $\SI{1.9e9}{}$. As to be expected, this condition number matches the corresponding number of the frequency-domain formulation for the same test problem when using $\omega=0$, see \cite{Herles_2025aa}. Again, the original system matrix is singular.

As a last example, we want to show that this time-domain formulation can be effectively used for nonlinear problems. We consider again the planar coil from section \ref{sec:numPlanarCoil}. This time we apply the nonlinear temperature dependent material law
\begin{equation}
    \conductivity(T) = \frac{\conductivity_0}{1 + \alpha (T-T_0)}, \label{eq:nonlinearConductivity}
\end{equation}
for the conductivity of copper with reference conductivity $\conductivity_0 = \SI{6e7}{\frac{S}{m}}$, temperature coefficient $\alpha = \SI{3.93e-3}{\frac{1}{K}}$ and reference temperature $T_0 = \SI{20}{^\circ C}$. 
Heat conduction in the copper coil is assumed to occur instantaneously. Thus, the temperature can be obtained by solving a lumped heat equation \cite{incropera2002fundamentals}.
Furthermore, a simplified coupling is realized by computing the electrical losses using
$$P_\text{eqs} = \int_{\Omega_\mathrm{C}}\conductivity \GRAD \phifield \cdot \GRAD \phifield \,\diffV\;,$$
i.e., neglecting the influence of the vector potential $\Afield$. This is possible in applications where the current flow is already determined with sufficient accuracy by the first step. This would be different, for example, in cases in which conductors are magnetically coupled.
The resulting temperature can then be used to update the conductivity according to \eqref{eq:nonlinearConductivity} which in turn can be used for computation of the second step. This form of coupling enables us to resolve the nonlinearity of the electrically dependent materials in the first cheap step while still being able to compute the inductive effects with the second step. We increase the voltage amplitude $\hat{u} = \SI{50}{V}$ and simulate the problem for $\mathcal{I} = \left( 0, \SI{6.67}{ms}\right)$ using $\deltaT = \SI{0.134}{ms}$. 

\begin{figure}
\centering
\pgfplotstableread[col sep=comma]{data/nonlinearSim.csv}\data

\begin{tikzpicture}
\begin{groupplot}[
  group style={
    group size=2 by 2,
    horizontal sep=1.5cm,
    vertical sep=2cm,
  },
  tudalineplot,
  height=10em,
  width=0.25\textwidth,
  xlabel={\scriptsize$t$ in s},
  xticklabel style={font=\scriptsize},
  yticklabel style={font=\scriptsize},
]

\nextgroupplot[
  ylabel={\scriptsize $V_1$ in V},
  ymin=-55, ymax=55,
  xmin=0, xmax=0.007,
]

\addplot+[mark size=.7mm,mark=*,mark options={fill=TUDa-3c,scale=0.3}, 
  TUDa-3c, line width=1pt] table[x index = 0, y index = 1] {\data};

\nextgroupplot[
  ylabel={\scriptsize $T$ in $^\circ$C},
  ymin=0, ymax=800,
  xmin=0, xmax=0.007,
]
\addplot+[mark size=.7mm,mark=*,mark options={fill=TUDa-3c,scale=0.3}, 
  TUDa-3c, line width=1pt] table[x index = 0, y index = 2] {\data};

  \nextgroupplot[
  ylabel={\scriptsize $P$ in W},
  ymin=0, ymax=2.5e6,
  xmin=0, xmax=0.007,
]
\addplot+[mark size=.7mm,mark=*,mark options={fill=TUDa-3c,scale=0.3}, 
  TUDa-3c, line width=1pt] table[x index = 0, y index = 3] {\data};

  \nextgroupplot[
  ylabel={\scriptsize $\conductivity$ in $\SI{}{\frac{S}{m}}$},
  ymin=0, ymax=6e7,
  xmin=0, xmax=0.007,
]
\addplot+[mark size=.7mm,mark=*,mark options={fill=TUDa-3c,scale=0.3}, 
  TUDa-3c, line width=1pt] table[x index = 0, y index = 4] {\data};
\end{groupplot}
\node[anchor=north] at ($(group c1r1.south)+(0,-1cm)$) {\small (a) Voltage Excitation};
\node[anchor=north] at ($(group c2r1.south)+(0,-1cm)$) {\small (b) Temperature};
\node[anchor=north] at ($(group c1r2.south)+(0,-1cm)$) {\small (c) Electric Loss Power};
\node[anchor=north] at ($(group c2r2.south)+(0,-1cm)$) {\small (d) Conductivity};
\end{tikzpicture}

\caption{Voltage, temperature, electric loss power and conductivity over time for planar coil example with nonlinear conductivity.}
\label{fig:nonLinearSim}
\end{figure}
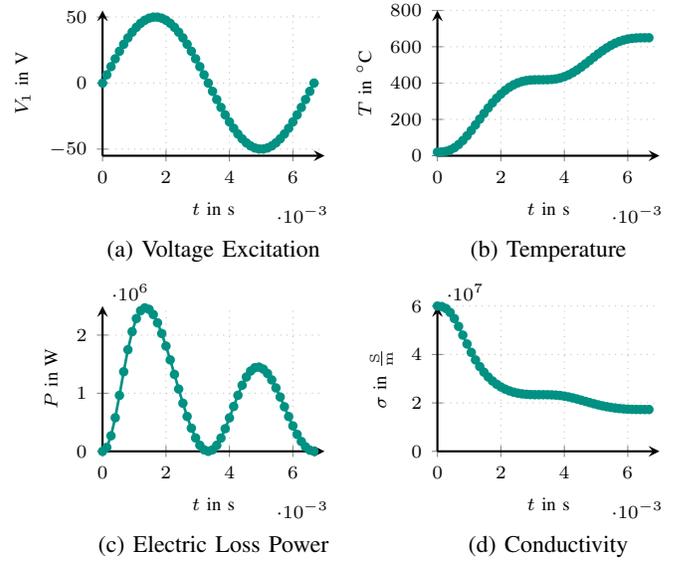

\autoref{fig:nonLinearSim} shows the sinusoidal voltage excitation as well as the temperature of the copper coil over the simulation time. One can see that the temperature rises significantly leading in turn to a nonlinear drop in conductivity. The resulting electric loss power decreases with time because of the lower current due to the temperature increase.

More generally, the option to consider nonlinear material dependencies in the cheap first step only is expected to be particularly interesting for surge-arresters, see e.g. \cite{Spack-Leigsnering_2019aa}, or for plasma applications, see e.g. \cite{busetto2024twostepmethodcouplingeddy}.

\section{Conclusion}\label{sec:conclusions}
We have presented a time-domain extension of the stabilized two-step formulation of Maxwell’s equations that enables numerically robust simulations in the low-frequency regime. The proposed method uses a generalized tree-cotree gauge to eliminate the kernel of the curl-curl operator, thus avoiding the instability commonly observed for large time steps. 

The first `electroquasistatic' step is discretized using the trapezoidal rule, while the second `full Maxwell' step is handled using the Newmark-beta method. Both schemes are second-order accurate and can be applied in parallel or subsequently, ensuring consistent time integration throughout the coupled system. 

Our numerical tests demonstrate that the approach yields accurate results even near the static limit, with significantly improved condition numbers compared to the original formulation. Furthermore, we present that nonlinearities in the conducting material can be treated very efficiently in cases where the first step already determines the current flow accurately. These results confirm that the stabilized time-domain two-step formulation is a promising framework for broadband and multiphysics electromagnetic simulations.

\section*{Acknowledgment}
The work is supported by the joint DFG/FWF Collaborative Research Centre CREATOR (DFG: Project-ID 492661287/TRR 361; FWF: 10.55776/F90) at TU Darmstadt, TU Graz and JKU Linz.

\bibliographystyle{IEEEtran}
\bibliography{bibtex}

\end{document}